\documentclass[12pt]{article}  % Comment this line out if you need a4paper

\usepackage{amsmath}
\usepackage{amssymb}
\usepackage{cite}
\usepackage{paralist}
\usepackage[usenames]{color}
\usepackage{graphicx}          % Include this line if your
\usepackage[dvips]{epsfig}    % or this line, depending on which
\newtheorem{definition}{Definition}
\newtheorem{theorem}{Theorem}
\newtheorem{corollary}{Corollary}

\newtheorem{assumption}{Assumption}

\title{On exponential stabilization of nonholonomic systems with time-varying drift\thanks{
This work is supported in part by the German Research Foundation (project GR 5293/1-1), the State Fund for Fundamental Research of Ukraine (project F78/206-2018), and NAS of Ukraine (budget program KPKBK 6541230).\newline
$^{1}$Institute of Mathematics, Julius Maximilian University of Würzburg,   Germany
        {\tt\small viktoriia.grushkovska@mathematik.uni-wuerzburg.de}
        \newline
$^{2}$Max Planck Institute for Dynamics of Complex Technical Systems, Magdeburg, Germany
{\tt\small zuyev@mpi-magdeburg.mpg.de}\newline
$^{3}$Institute of Applied Mathematics and Mechanics, National Academy of Sciences of Ukraine
}
}

\author{Victoria Grushkovskaya$^{1,3}$ and Alexander Zuyev$^{2,3}$}
\date{}

\begin{document}

\maketitle
\thispagestyle{empty}
%\pagestyle{empty}

%%%%%%%%%%%%%%%%%%%%%%%%%%%%%%%%%%%%%%%%%%%%%%%%%%%%%%%%%%%%%%%%%%%%%%%%%%%%%%%%
\begin{abstract}
A class of nonlinear control-affine systems with bounded time-varying drift is considered.
It is assumed that the control vector fields together with their iterated Lie brackets satisfy H\"ormander's condition in a neighborhood of the origin.
Then the problem of exponential stabilization is treated by exploiting periodic time-varying feedback controls.
An explicit parametrization of such controllers is proposed under a suitable non-resonance assumption.
It is shown that these controllers ensure the exponential stability of the closed-loop system provided that the period is small enough.
The proposed control design methodology is applied for the stabilization of an underwater vehicle model and a front-wheel drive car.
\end{abstract}

%%%%%%%%%%%%%%%%%%%%%%%%%%%%%%%%%%%%%%%%%%%%%%%%%%%%%%%%%%%%%%%%%%%%%%%%%%%%%%%%
\section{Introduction}

The paper focuses on the stabilization problem for a class of nonholonomic systems in the control-affine form.
As the number of control inputs in such systems can be significantly smaller than the dimension of the state vector, this causes certain challenges in control design.
 There exists a number of approaches which allow to stabilize \emph{control-linear} nonholonomic systems  (see, e.g.,~\cite{Cor92,Ast94,Mor99,ZuSIAM}, and references therein). However, the stabilization problem becomes even more complicated for \emph{control-affine} systems with  unstable drift terms.
Controllability properties and motion planning problems of control-affine systems were discussed, e.g., in~\cite{Luca95,Godv99,Pom99,Ag12,Jean15,ZG17}. While rather general results have been obtained for motion planning problems, stabilization of nonholonomic systems with drift is mainly studied for specific classes of systems~(see, e.g., \cite{MCl98,Rey99,Bu00,Thi00,Wang04,YaYa10,Gao11,Zhao13}, and~\cite{Kolman95,Mich03} for a survey).
A more general class of control-affine systems was considered in~\cite{Her80,Mich03}, where stabilizing controllers have been proposed under the assumption that the system is strongly controllable and can be approximated by a system with  nilpotent Lie algebra, and that the drift term vanishes at the origin.
\\
In this paper, we propose a class of control functions that  stabilize the origin of an underactuated control-affine system with time-varying drift term. In general, we do not assume that the drift vanishes at the origin, which leads to the practical asymptotic stability of the corresponding closed-loop system. For a special class of drift terms vanishing at the origin, we show that the trajectories of the system exponentially tend to zero.
We also do not involve the drift vector field in the controllability rank condition. %, rather considering it as a disturbance or a part of system's dynamics.
 In Section~2, we formulate the problem statement and present a novel stabilizability result as an the extension of the control design approach from~(\cite{ZuSIAM,GZ18}). Section~3 contains the proofs. Several examples are presented in Section~4.

\section{Main results}
\subsection{Problem statement}
Consider a system
\begin{equation}\label{sys1}
  \dot x= g(t,x)+\sum_{i=1}^mf_i(x)u_i,\quad x\in D\subset \mathbb R^n,
\end{equation}
where $x=(x_1,\dots,x_n)^\top$ is the state, $u=(u_1,\dots,u_m)^\top\in\mathbb R^m$ is the control, $f_i\in C^3(D;\mathbb R^n)$ describe the system dynamics,  and $g:\mathbb R\times D\to\mathbb R^n$ is the drift term related to the system dynamics or to disturbances.
In this paper, we propose a family of control laws for stabilizing the origin of system~\eqref{sys1} under the assumption that the vector fields $f_i$ together with their first- and second-order Lie brackets span the whole $n$-dimensional space, and the drift $g$ satisfies certain boundedness assumptions.
\begin{assumption}[Rank condition]
Let
$$S_1\subseteq \{1,2,...,m\},\;  S_2\subseteq \{1,2,...,m\}^2,\; S_3 \subseteq\{1,2,...,m\}^3$$ be  sets of indices  such that $|S_1|+|S_2|+|S_3|=n$ and, for each $x\in D$,
\begin{equation}\label{rank}
\begin{aligned}
 {\rm span}\big\{f_{i}(x), &[f_{j_1},f_{j_2}](x), {\left[[f_{\ell_1},f_{\ell_2}],f_{\ell_3}\right](x)} \,|\,\\
 & i{\in}S_1,
(j_1,j_2){\in} S_2, (\ell_1,\ell_2,\ell_3){\in} S_3\big\}{=}\mathbb{R}^n.
\end{aligned}
\end{equation}
\end{assumption}
\begin{assumption}[Boundedness of the drift]
For each  compact set $\xi\subseteq D$, there exists a $\tau>0$ and $\mathcal M_g\ge 0$ such that, for any $t_0\ge0$,
$
\|g(t,x)\|\le \mathcal M_g\text{ for all }t\in[t_0,t_0+\tau],\,x\in \xi.
$
\end{assumption}
To stabilize system~\eqref{sys1} at $x^*=0$, we adopt the control design approach previously  proposed for the case $g(t,x)=0$ in~\cite{ZGB16,GZ18}. Note that the presence of non-zero drift may  affect significantly the system behavior and complicates  the stabilization problem. Therefore, the results of the above mentioned papers cannot be directly applied, and  more sophisticated analysis is required.
\subsection{Notations and definitions}
\begin{definition}
We say that there is a \emph{resonance of order $N\in\mathbb N$} between the pairwise distinct  numbers $k_1,\dots,k_n$, if there exist relatively prime integers  $c_1,{\dots},c_n$  such that $|c_1|+...+|c_n|=N$ and
$c_1k_1+...+c_nk_n=0$.
\end{definition}
Similarly to the approaches of~\cite{Clar97,ZuSIAM}, we will exploit the sampling concept.
For a given $\varepsilon>0$, define a partition $\pi_\varepsilon$ of $[0,+\infty)$ into the intervals
$
[t_j,t_{j+1})$, $t_j=\varepsilon j$, $j=0,1,2,\dots $ .
\begin{definition}
\emph{Given a feedback $u=h(t,x)$, $h:[0,+\infty)\times D\to\mathbb R^m$, $\varepsilon>0$, and $x^0\in D$, a $\pi_\varepsilon$-solution of~\eqref{sys1} corresponding to $x^0$ and $h(t,x)$ is an absolutely continuous function  $x(t)\in D$, defined for $t\in[0,+\infty)$, such that  $x(0)=x^0$ and
$
\dot x(t)=f\big(x(t), h(t,x(t_j))\big)$, $t\in[t_j,t_{j+1}),
$
for each j=0,1,2,\dots.}
\end{definition}
  For  $f,g:\mathbb R^n\to\mathbb R^n $, $x\in\mathbb R^n$, the directional derivative is denoted as
 {$ L_gf(x)=\lim\limits_{s\to0}\tfrac{f(x+sg(x))-f(x)}{s}$}, and  $[f,g](x)= L_fg(x)- L_gf(x)$ stands for the Lie bracket.
 Throughout this paper, $\|a\|$ denotes the Euclidean norm of a vector $a\in\mathbb R^n$, and the norm of an $n\times n$-matrix $\cal F$ is defined as $\|{\cal F}\|=\sup_{\|y\|=1}\|{\cal F}y\|$.
\subsection{Control functions}
 Given positive real numbers $\varepsilon$ and $\gamma$, we define the control functions $u_k$, $k=1,\dots,m,$ as
\begin{align}
u_k= h&^\varepsilon_k(t,x)=\sum_{i_1\in S_1} a_{i_1}(x)\phi^{(k,\varepsilon)}_{i_1}(t)\nonumber \\
 &+\varepsilon^{-\tfrac{1}{2}}\sum_{(j_1,j_2)\in S_2} \sqrt{|a_{j_1j_2}(x)|}\phi^{(k,\varepsilon)}_{j_1j_2}(t,x)\label{cont}\\
  &+\varepsilon^{-\tfrac{2}{3}}\sum_{(\ell_1,\ell_2,\ell_3)\in S_3} \sqrt[3]{a_{\ell_1\ell_2\ell_3}(x)}\phi^{(k,\varepsilon)}_{\ell_1\ell_2\ell_3}(t),\nonumber
\end{align}
where the  state-dependent vector function
\begin{align*}
a(x)&=\big(a_{i_1}(x)\big|_{i_1\in S_1}, a_{j_1j_2}(x)\big|_{(j_1,j_2)\in S_2}, \nonumber\\
&\qquad\qquad\qquad a_{\ell_1\ell_2\ell_3}(x)\big|_{(\ell_1,\ell_2,\ell_3)\in S_3}\big)^\top\in\mathbb R^n\nonumber
\end{align*}
is chosen as
\begin{align}
a(x)=- \gamma \mathcal F^{-1}(x) x\label{a}
\end{align}
with some control gain $\gamma>0$,  and $\phi^{(k,\varepsilon)}_{i_1}(t)=\delta_{k i_1},$
\begin{align}
 \phi^{(k,\varepsilon)}_{j_1j_2}(t,x)&{=}2\sqrt{\pi \kappa_{j_1j_2}}\Big(\delta_{kj_1}{\rm sign}(a_{j_1,j_2}(x))\cos{\frac{2\pi \kappa_{j_1j_2}}{\varepsilon}}t\nonumber\\
 &{+}\delta_{kj_2}\sin{\frac{2\pi \kappa_{j_1j_2}}{\varepsilon}}t\Big), \nonumber\\
\phi^{(k,\varepsilon)}_{\ell_1\ell_2\ell_3}(t)&{=}2\sqrt[3]{2\pi^2 \kappa_{3\ell_1\ell_2\ell_3}\kappa_{4\ell_1\ell_2\ell_3}}\Big(\delta_{kl3_1}\cos{\frac{2\pi \kappa_{1\ell_1\ell_2\ell_3}t}{\varepsilon}}\nonumber\\
&{+}\delta_{k\ell_2}\sin{\frac{2\pi \kappa_{2\ell_1\ell_2\ell_3}t}{\varepsilon}} \nonumber\\
&{+}\delta_{k\ell_3}{\cos{\frac{2\pi \kappa_{1\ell_1\ell_2\ell_3}t}{\varepsilon}}\sin{\frac{2\pi \kappa_{2\ell_1\ell_2\ell_3}t}{\varepsilon}}}\Big).\label{phi}
\end{align}
Here $\delta_{ki}$ is the Kronecker delta, and the integer parameters
$\kappa_{j_1j_2}$, $\kappa_{1\ell_1\ell_2\ell_3}$, $\kappa_{2\ell_1\ell_2\ell_3}$ are specified according to the following assumption.
\begin{assumption}[Absence of resonances]~\label{ControlDesign}
\emph{The positive integer numbers  $\kappa_{j_1j_2}$, $\kappa_{1\ell_1\ell_2\ell_3}$, $\kappa_{2\ell_1\ell_2\ell_3}$, $\kappa_{3\ell_1\ell_2\ell_3}=\kappa_{1\ell_1\ell_2\ell_3}+\kappa_{2\ell_1\ell_2\ell_3}$, and $\kappa_{4\ell_1\ell_2\ell_3}=\kappa_{2\ell_1\ell_2\ell_3}-\kappa_{1\ell_1\ell_2\ell_3}$ are pairwise distinct, and there are no  third-order resonances between $\kappa_{s\ell_1\ell_2\ell_3}$ ($s=1,\dots,4$), except those imposed by the definition of $\kappa_{3\ell_1\ell_2\ell_3},\kappa_{4\ell_1\ell_2\ell_3}$.}
  \end{assumption}

\subsection{Stabilization of system~\eqref{sys1}}
Consider the matrix
\begin{equation}\label{fmatrix}
\begin{aligned}
 \mathcal F(x)= \Big(\big(f_{i}(x)&\big)_{j_1\in S_1}\ \big([f_{j_1},f_{j_2}](x)\big)_{(j_1,j_2)\in S_2} \\
 &\big(\left[[f_{\ell_1},f_{\ell_2}],f_{\ell_3}\right](x)\big)_{(\ell_1,\ell_2,\ell_3)\in S_3}\Big),
\end{aligned}
\end{equation}
which is nonsingular in $D$ provided that condition~\eqref{rank} holds.
The main result of this paper is the following theorem.
\begin{theorem}\label{thm_practical}
  \emph{Let $D\subseteq \mathbb R^n$, $f_i\in C^3(D;\mathbb R^n)$, $i=1,\dots,m$. Suppose that  Assumptions~1--2 hold in $D$ and  there exists an  $\alpha>0$ such that
 $ \|\mathcal F^{-1}(x)\|\le \alpha \text{ for all }x\in D,$
where the matrix $\mathcal F(x)$ is given by~\eqref{fmatrix}.}
\\
\emph{
If the functions $u_k=h^\varepsilon_k(t,x)$, $k=1,\dots,m$, are defined as in~\eqref{cont}--\eqref{phi} with the parameters satisfying Assumption~\ref{ControlDesign}, then
for any $\delta,\rho>0$  there exist  $\gamma,\bar\varepsilon>0$  such that, for any $\varepsilon\in(0,\bar\varepsilon]$, the $\pi_\varepsilon$-solution of system~\eqref{sys1} with the initial data $x(0)=x^0\in B_\delta(0)$ is well-defined on $t\in [0,+\infty)$ and
$$
\begin{aligned}
 \|x(t)\|\le \|x^0\|e^{-\lambda t} +\rho&\text{ for all } t\in[t_0,t_1),\\
 \text{and }\|x(t)\|\le \rho&\text{ for all } t\in[t_1,\infty),\\
\end{aligned}
 $$
 with some $\lambda,t_1>0$.}
\end{theorem}
The proof  is given in  Section~\ref{proofs}.1. Note that the proof provides a constructive procedure for choosing $\gamma,\lambda$ and $\bar\varepsilon$. Theorem~\ref{thm_practical} gives the \emph{practical exponential stability} conditions of the point $x=0$.
Obviously, to stabilize system~(1) in the practical sense at any other point $x^*\in D$,  one can take $a(x)=- \gamma \mathcal F^{-1}(x) (x-x^*)$.
Under some  stronger assumptions on $g(t,x)$, even local exponential  stability can be achieved, as stated in the following corollaries.
\begin{corollary}
\emph{Let $D\subseteq \mathbb R^n$, $f_i\in C^3(D;\mathbb R^n)$, $i=1,\dots,m$. Assume that  Assumption~1 holds in $D$ and  there exists an  $\alpha>0$ such that
 $ \|\mathcal F^{-1}(x)\|\le \alpha \text{ for all }x\in D,$
where the matrix $\mathcal F(x)$ is given by~\eqref{fmatrix}. Assume also that there are $\mathcal M_g,\mathcal L_g\ge 0$ and $\delta_0>0$ such that
$$g(t,x)\le \mathcal M_g\|x\|^3,\,\Big\|g(t,x)-g(t,y)\Big\|\le \mathcal L_g\|x-y\|,$$
for all $t\ge0,x,y\in B_{\delta_0}(0).$}
\emph{If the functions $u_k=h^\varepsilon_k(t,x)$, $k=1,\dots,m$, are defined as in~\eqref{cont}--\eqref{phi} with the parameters satisfying Assumption~\ref{ControlDesign}, then
for any $\delta>0$  there exist  $\gamma,\bar\varepsilon>0$  such that, for any $\varepsilon\in(0,\bar\varepsilon]$, the $\pi_\varepsilon$-solution of system~\eqref{sys1} with the initial data $x(0)=x^0\in B_\delta(0)$ is well-defined on $t\in [0,+\infty)$ and
$$
 \|x(t)\|=O(e^{-\lambda t})\text{ as } t\to\infty,\text{ with some }\lambda>0.
 $$
}
\end{corollary}
The proof of Corollary~1   is  in  Section~\ref{proofs}.2.
\section{Proofs of the main results}~\label{proofs}
\subsection{Proof of Theorem~\ref{thm_practical}}
For any $x^0\in D$,  let $\rho,\delta,\delta'>0$ be such that $\overline{B_{\rho}(0)}\subset\overline{B_\delta(0)}\subset \overline{B_{\delta'}(0)}\subset D,$ and
$
\mathcal M_f=\sup_{x\in \overline{B_{\delta'}(0)}}\|f_i(x)\|$, $\mathcal M_g=\sup_{x\in \overline{B_{\delta'}(0)},t\in[0,\tau]}\|g(t,x)\|.
$
Let $\varepsilon_0=\min\Big\{\tau,\frac{1}{\gamma}\Big\}$ and
$
U^\varepsilon(x^0)=\max\limits_{0\le t\le \varepsilon}\sum_{i=1}^m |h_i^{\varepsilon}(t,x^0)|.
$ Here we assume that $\gamma>0$ is fixed, since, as it will be shown later,  $\gamma$ can be defined independently on $\varepsilon$.
 From~(\cite{GZ18}), for every $\varepsilon\in(0,\varepsilon_0)$,
\begin{equation}\label{u_est}
\begin{aligned}
U^\varepsilon(x^0)&\le c_{u1}\gamma\|x^0\|+c_{u2}\sqrt{\frac{\gamma}{\varepsilon}\|x^0\|}+\sqrt[3]{\frac{\gamma}{\varepsilon^2}\|x^0\|},\\
U^\varepsilon\varepsilon(x^0)&\le c_u\sqrt[3]{\varepsilon\gamma\|x^0\|},
\end{aligned}
\end{equation}
where
$$c_{u1}= \alpha\sqrt{|S_1|},c_{u2}=4\sqrt{{\pi\alpha}}\Big(\sum\limits_{(j_1,j_2)\in S_2}{\kappa_{j_1j_2}}^{2/3}\Big)^{3/4},$$
 $$c_{u3}{=}6\sqrt[3]{2\pi^2\alpha}\Big(\sum\limits_{(\ell_1,\ell_2,\ell_3)\in S_3}|\kappa_{2\ell_1\ell_2\ell_3}^2{-}\kappa_{1\ell_1\ell_2\ell_3}^2|^{2/5}\Big)^{5/6},$$
 and% $c_u=c_{u1}+c_{u2}+c_{u3}$.
 $$
 c_u=c_{u1}\|x^0\|^{2/3}+c_{u2}\|x^0\|^{1/6}+c_{u3}.
 $$
The integral representation
$$
x(t)-x^0=\int_0^t\Big(\sum_{i=1}^mf_i(x(s))h_i^\varepsilon(s,x^0)+g(s,x(s))\Big)ds
$$
yields that, for any $x^0\in B_{\delta}(0)$, $\varepsilon\in(0,\varepsilon_0)$,
$$
\|x(t)-x^0\|\le \mathcal M_fc_u\sqrt[3]{\varepsilon\gamma\|x^0\|}+\varepsilon \mathcal M_g\text{ for all }t\in[0,\varepsilon].
$$
For $d=\min\{\delta'-\delta,\tfrac{1}{2}\rho\}>0$, let $\varepsilon_1$ be the smallest positive root of the equation
$$
\mathcal M_fc_u\sqrt[3]{\varepsilon\gamma\delta'}+\varepsilon \mathcal M_g=d.
$$
Then for any $\varepsilon\in(0,\min\{\varepsilon_0,\varepsilon_1\})$, the solutions of~\eqref{sys1}, \eqref{cont} with $x(0)\in B_\delta$ are well defined in $D$ ($\|x(t)\|\le \delta'$) for $t\in[0,\varepsilon]$, and
\begin{equation}\label{xless}
\text{if }\|x^0\|\le \frac{\rho}{2} \text{ then }\|x(t)\|\le \rho\text{ for all }t\in[0,\varepsilon].
\end{equation}
Then we use the Chen--Fliess series to represent the $\pi_\varepsilon$-solution of system~\eqref{sys1} at time $\varepsilon$, taking into account the drift term $g(t,x)$ and  formula~\eqref{a}:
\begin{equation}\label{volt_drift}
\begin{aligned}
x(\varepsilon)=&x^0+\varepsilon\sum_{j_1\in S_1}f_{j_1}(x^0)a_{j_1}(x^0)\\
&+\varepsilon\sum_{(j_1,j_2)\in S_2}[f_{j_1},f_{j_2}](x^0)a_{j_1j_2}(x^0)\\
&+ \varepsilon\sum_{(\ell_1,\ell_2,\ell_3)\in S_3}[[f_{\ell_1},f_{\ell_2}],f_{\ell_3}](x^0)a_{\ell_1\ell_2\ell_3}(x^0)\\
&+\Omega(a,\varepsilon)+r_f(\varepsilon)+r_g(\varepsilon)+\int_0^\varepsilon g(s,x(s))ds\\
= &x^0-\gamma\varepsilon x^0+\int_0^\varepsilon g(s,x(s))ds+\Omega(a,\varepsilon)\\
&+r_g(\varepsilon)+r_a(\varepsilon),
\end{aligned}
\end{equation}
$$
\begin{aligned}
r_g&(\varepsilon)=\int\limits_0^\varepsilon\int\limits_0^{s_1}\sum_{j_1=1}^m L_gf_{j_1}(x(s))h_{j_1}^\varepsilon(s_1,x^0){\rm d}s_2{\rm d}s_1\\
&+\int\limits_0^\varepsilon\int\limits_0^{s_1}\int\limits_0^{s_2}\sum_{j_1,j_2=1}^m L_g L_{f_{j_2}}f_{j_1}(x(p))\\
&\quad\quad\quad\times h_{j_2}^\varepsilon(s_2,x^0)h_{j_1}^\varepsilon(s_1,x^0){\rm d}s_3{\rm d}s_2{\rm d}s_1\\
&+\sum\limits_{\hspace{-1em}j_1,j_2,j_3=1}^m{\int\limits_0^\varepsilon}{\int\limits_0^{s_1}}{\int\limits_0^{s_2}}{\int\limits_0^{s_3}} L_g L_{f_{j_3}} L_{f_{j_2}}f_{j_1}(x(s_4))\\
 & \quad\quad\quad\times h_{j_3}^\varepsilon(s_3,x^0)  h_{j_2}^\varepsilon(s_2,x^0)h_{j_1}^\varepsilon(s_1,x^0){\rm d}s_4{\rm d}s_3{\rm d}s_2 {\rm d}s_1,\\
r_f& (\varepsilon)={\sum\limits_{\hspace{-1em}j_1,\dots,j_4=1}^m}{\int\limits_0^\varepsilon}{\int\limits_0^{s_1}}{\int\limits_0^{s_2}}{\int\limits_0^{s_3}} L_{f_{j_4}}L_{f_{j_3}} L_{f_{j_2}}  f_{j_1}(x(s_4))
 \\
 \times&
 h_{j_4}^\varepsilon(s_4,x^0)  h_{j_3}^\varepsilon(s_3,x^0)  h_{j_2}^\varepsilon(s_2,x^0)  h_{j_1}^\varepsilon(s_1,x^0){\rm d}s_4{\rm d}s_3{\rm d}s_2{\rm d}s_1.
\end{aligned}
$$
We omit the explicit expression for $\Omega(a,\varepsilon)$ due to the space limits. Similarly to~(\cite{GZ18,GZE18,ZG19}), it can be shown that there exist $c_\Omega,c_g,c_f\ge0$ such that, for any $x^0\in \overline{B_\delta(0)}$,
$$
\begin{aligned}
\|\Omega(a,\varepsilon)\|&\le c_\Omega(\varepsilon\|x^0\|)^{7/6},\\
\|r_g(\varepsilon)\|&\le c_gM_g\varepsilon^{4/3}\|x^0\|^{1/3},\ \|r_f(\varepsilon)\|\le c_f(\varepsilon\|x^0\|)^{4/3}.\\
\end{aligned}
$$
Applying these estimates to~\eqref{volt_drift}, we conclude that
\begin{equation}\label{x_est}
\|x(\varepsilon)\|\le (1-\gamma\varepsilon)\|x^0\|+\sigma(\varepsilon) \varepsilon^{7/6}\|x^0\|^{1/3}+\mathcal M_g\varepsilon,
\end{equation}
where $\sigma(\varepsilon)=c_\Omega\delta^{5/6}+\varepsilon^{1/6}\big(c_gM_g+c_f\delta\big)$.
Assume $x^0\in \overline{B_\delta(0)}\setminus B_{\rho/2}(0)$. Then the latter inequality can be rewritten as
$$
\begin{aligned}
\|x(\varepsilon)\|&\le (1-\gamma\varepsilon)\|x^0\|+\sigma \varepsilon^{7/6}\Big(\frac{2}{\rho}\Big)^{2/3}\|x^0\|+\frac{2 \mathcal M_g}{\rho}\varepsilon\|x^0\|\\
&=(1-\varepsilon\lambda_1)\|x^0\|,
\end{aligned}
$$
where $\lambda_1= \gamma-\frac{2 \mathcal M_g}{\rho}-\sigma(\varepsilon) \varepsilon^{1/6}\Big(\frac{2}{\rho}\Big)^{2/3}$. Taking $\gamma>\frac{2\mathcal M_g}{\rho}$, we ensure that  there exists a $\lambda_2>0$ such that
$\gamma-\frac{2 \mathcal M_g}{\rho}>\lambda_2$.  For any $\lambda\in(0,\lambda_2)$, let $\varepsilon_2=\min\Big\{\frac{1}{\lambda},\hat\varepsilon\Big\}$, where $\hat\varepsilon$ is the smallest positive root of the equation
$$
\sigma(\varepsilon) \varepsilon^{1/6}\Big(\frac{2}{\rho}\Big)^{2/3}=\lambda_2-\lambda.
$$
Then, for any $\varepsilon\in(0,\bar\varepsilon=\min\{\varepsilon_0,\varepsilon_1,\varepsilon_2\})$, if $\|x^0\|>\frac{\rho}{2}$, then
$$
\|x(\varepsilon)\|\le (1-\varepsilon\lambda)\|x^0\|.
$$
Since $x^0\in B_\delta(0)$ then  $x(\varepsilon)\in B_\delta(0)$, and we  repeat the above argumentation for  the solutions $x(t)$  of system~\eqref{sys1}, \eqref{cont} with the initial conditions $x(\varepsilon)\in B_\delta(0)$. Thus, we  conclude that there exists an $N\in\mathbb N\cup \{0\}$ such that
$$
\begin{aligned}
\|x(j\varepsilon)\|\ge\frac{\rho}{2}\text{ for all }j=0,\dots,N-1,\ \|x(N\varepsilon)\|\le\frac{\rho}{2},
\end{aligned}
$$
which implies that the solutions $x(t)$  of system~\eqref{sys1}, \eqref{cont} with the initial conditions $x(0)=x^0\in B_\delta(0)$ are well defined for all $t\in[0,(N+1)\varepsilon]$, and
$$
\|x(j\varepsilon)\|\le \|x^0\|e^{-\lambda j\varepsilon} \text{ for all }j=0,\dots,N.
$$
Furthermore, $\|x\big((N+1)\varepsilon\big)\|\le \rho$ from~\eqref{xless}. If  $\|x\big((N+1)\varepsilon\big)\|\ge\frac{\rho}{2}$, we apply again the same reasoning and obtain $\big\|x\big((N+2)\varepsilon\big)\big\|\le \|x\big((N+1)\varepsilon\big)\|$. Otherwise,~\eqref{xless} implies $\|x\big((N+2)\varepsilon\big)\|\le \rho$. Thus, for any $\varepsilon\in(0,\bar\varepsilon)$, the solutions of system~\eqref{sys1}, \eqref{cont} with the initial conditions $x(0)=x^0\in B_\delta(0)$ satisfy the following properties:
$$
\|x(0)\|\le \delta\Rightarrow\|x(t)\|\le \|x^0\|e^{-\lambda t} +\frac{\rho}{2}\text{ for all }t\ge 0,
$$
and there exists a $t_1>0$ such that
$\|x(t)\|\le \rho \text{ for }t\ge t_1.
$
\subsection{Proof of Corollary~1}
As it follows from  Theorem~\ref{thm_practical} and its proof, for any $\delta,\delta_0>0$, there exists an $\bar\varepsilon_1>0$ such that, for any $\varepsilon\in(0,\bar\varepsilon_1]$, the $\pi_\varepsilon$-solution of system~\eqref{sys1} with the initial data $x(0)=x^0\in B_\delta(0)$ is well-defined on $t\in [0,+\infty)$ and
\begin{equation}\label{exp1}
\begin{aligned}
 \|x(t)\|\le \|x^0\|e^{-\lambda_1 t} +\delta_0&\text{ for all } t\in[t_0,t_1),\\
 \text{and }\|x(t)\|\le \delta_0&\text{ for all } t\in[t_1,\infty),\\
\end{aligned}
\end{equation}
 with some $\lambda_1,t_1>0$.
The proof is similar to the proof of Theorem~\ref{thm_practical} with $D=B_{\delta_0}(0)$, so we just briefly describe the main differences. Let us analyze the behavior of solutions of system~\eqref{sys1} in $B_{\delta_0}(0)$.

Let $\tilde x^0\in B_{\delta_0}(0)$. Using the integral representation of $x(t)$,  the Gr\"onwall–-Bellman inequality, estimate~\eqref{u_est}, and the assumptions on $g(t,x)$, we conclude that
%Extending~\cite[Lemma~4]{ZuSIAM} to system~\eqref{sys1}, we obtain  the following estimate:
%$$
%\|x(t){-}\tilde x^0\|{\le}\frac{\big(\mathcal M_fU^\varepsilon(\tilde x^0){+}\mathcal M_g\|\tilde x^0\|^3\big)}{\mathcal L_fU^\varepsilon(\tilde x^0){+}\mathcal L_g}\Big(e^{\big(\mathcal L_fU^\varepsilon(\tilde x^0){+}\mathcal L_g\big)t}{-}1\Big),
%$$
\begin{equation}\label{xx0}
  \|x(t)-\tilde x^0\|\le c_x\sqrt[3]{\varepsilon\|\tilde x^0\|}\text{ for all }t\in[0,\varepsilon],
\end{equation}
where $c_x=\big(\mathcal M_fc_u\sqrt[3]\gamma+\mathcal M_g\delta_0^2(\varepsilon\delta_0)^{2/3}\big)e^{\mathcal L_fc_u\sqrt[3]{\varepsilon\gamma\delta_0}+\mathcal L_g}$, and $\mathcal L_f$ is such that $\Big\|f(x)-f(y)\Big\|\le \mathcal L_f\|x-y\|$ for all $x,y\in B_{\delta'}(0).$
%It can be shown that there exists an $\tilde \varepsilon_0>0$ such that, for any $\varepsilon\in(0,\tilde\varepsilon_0)$, the solutions of~\eqref{sys1}, \eqref{cont} with $x(0)\in B_\delta$ are well defined in $D$, and
Furthermore,
\begin{equation}\label{g}
\begin{aligned}
  \|g(t,&x(t))\|\le \mathcal M_g\|x(t)\|^3\le \mathcal M_g(\|\tilde x^0\|+\|x(t)-\tilde x^0\|)^3\\
  &\le \mathcal M_g\|\tilde x^0\|(\delta_0^2+c_x\sqrt[3]{\varepsilon})^3\text{ for all }t\in[0,\varepsilon].
\end{aligned}
\end{equation}
Then the term $r_g(\varepsilon)$    in~\eqref{volt_drift} can be estimated as
$
\|r_g(\varepsilon)\|\le \tilde c_g (\varepsilon\|\tilde x^0\|)^{4/3}
$
with some $\tilde c_g>0$. Consequently,
  the estimate~\eqref{x_est} can be written as
$$
\begin{aligned}
\|x(\varepsilon)\|\le& (1-\gamma\varepsilon)\|x^0\|+\tilde\sigma(\varepsilon) \varepsilon^{7/6}\|x^0\|^{7/6}\\
&+\varepsilon \mathcal M_g\|\tilde x^0\|(\delta_0^2+c_x\sqrt[3]{\varepsilon})^3\\
=&\Big(1-\varepsilon\big(\gamma-\mathcal M_g\delta_0^6-\varepsilon^{1/6}\sigma_1(\varepsilon)\big)\Big)\|x^0\|.
\end{aligned}
$$
Here $\tilde\sigma(\varepsilon)=c_\Omega+(\varepsilon\delta_0)^{1/6}\big(\tilde c_g+c_f\big)$, $\sigma_1(\varepsilon)=\tilde\sigma(\varepsilon)\delta_0^{1/6}$. Taking $\gamma>\mathcal M_g\delta_0^6$, $\lambda_2\in(0,\gamma-\mathcal M_g\delta_0)$, and $\tilde\varepsilon_1$ as the smallest positive root of the equation $\varepsilon^{1/6}\sigma_1(\varepsilon)=\tilde\lambda$, we obtain
$
\|x(\varepsilon)\|\le (1-\lambda_2\varepsilon)\|\tilde x^0\|.
$
Repeating the above argumentation for an arbitrary $\tilde x^0\in B_{\delta_0}(0)$, we conclude that
\begin{equation}\label{xjeps}
\|x(j\varepsilon)\|\le \|\tilde x^0\|e^{-\lambda_2 j\varepsilon}, \text{ for }j=0,1,2,\dots\, .
\end{equation}
For any $t\ge 0$ and $\varepsilon\in(0,\bar\varepsilon=\min\{\tilde\varepsilon_0,\tilde\varepsilon_1\})$, we have
$$
\begin{aligned}
\|x(t)\|&\le \Big\|x(t)-x\Big(\Big[\frac{t}{\varepsilon}\Big]\varepsilon\Big)\Big\|+\Big\|x\Big(\Big[\frac{t}{\varepsilon}\Big]\varepsilon\Big)\Big\|\\
&\le \sqrt[3]{\Big\|x\Big(\Big[\frac{t}{\varepsilon}\Big]\varepsilon\Big)\Big\|}\Big( c_x\sqrt[3]{\varepsilon}+\Big\|x\Big(\Big[\frac{t}{\varepsilon}\Big]\varepsilon\Big)\Big\|^{2/3}\Big).
\end{aligned}
$$
Using~\eqref{xjeps}, we obtain the following estimate:
\begin{equation}\label{exp2}
\|x(t)\|\le\mu_1 \sqrt[3]{\|\tilde x^0\|}e^{-\frac{\lambda_2}{3}t},
\end{equation}
with $\mu_1=e^{\lambda_2\varepsilon}\big(c_x\sqrt[3]{\varepsilon}+\delta_0^{2/3}\big)$.
Choosing $\bar\varepsilon=\min\{\bar\varepsilon_1,\bar\varepsilon_2\}$ and summarizing~\eqref{exp1} and~\eqref{exp2}, we conclude that, for any $\varepsilon\in(0,\bar\varepsilon]$, there exists a $ t_1> 0$
$$
\|x(t)\|\le\left\{
\begin{aligned}
  &\|x^0\|e^{-\lambda_1 t} +\delta_0\text{ for }t\in[0,t_1),\\
  &\mu_1 \sqrt[3]{\|x(t_1)\|}e^{-\frac{\lambda_2}{3}t}\text{ for }t\in[t_1,\infty),
\end{aligned}
\right.
$$
which proves the Corollary.
\section{Examples}
\subsection{Underwater vehicle with drift}
\begin{figure*}[tpt]
 \begin{minipage}{0.49\linewidth}
\begin{center}
\includegraphics[width=1\linewidth]{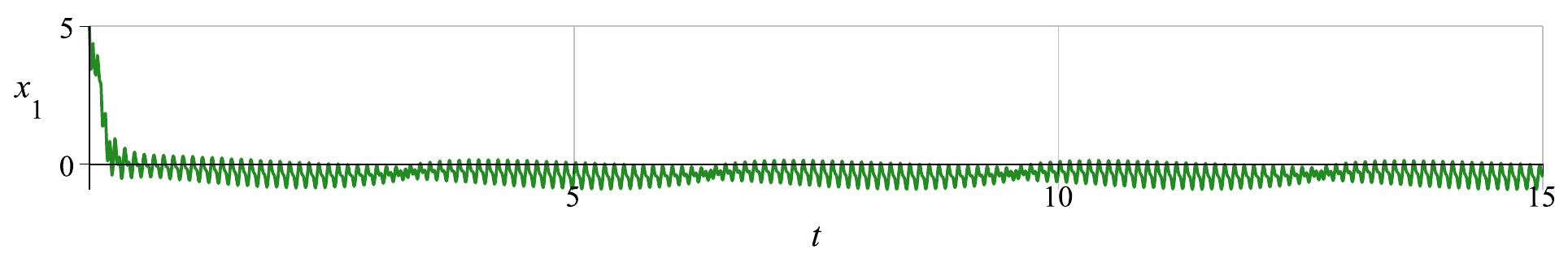}
\includegraphics[width=1\linewidth]{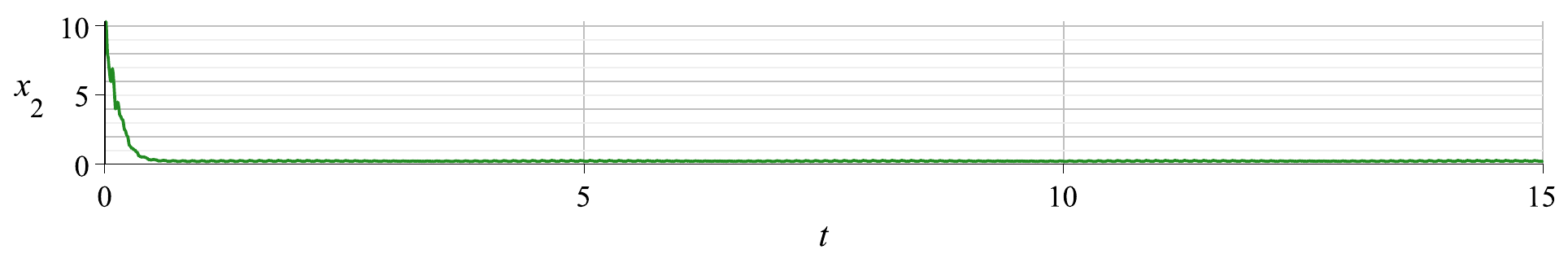}
\includegraphics[width=1\linewidth]{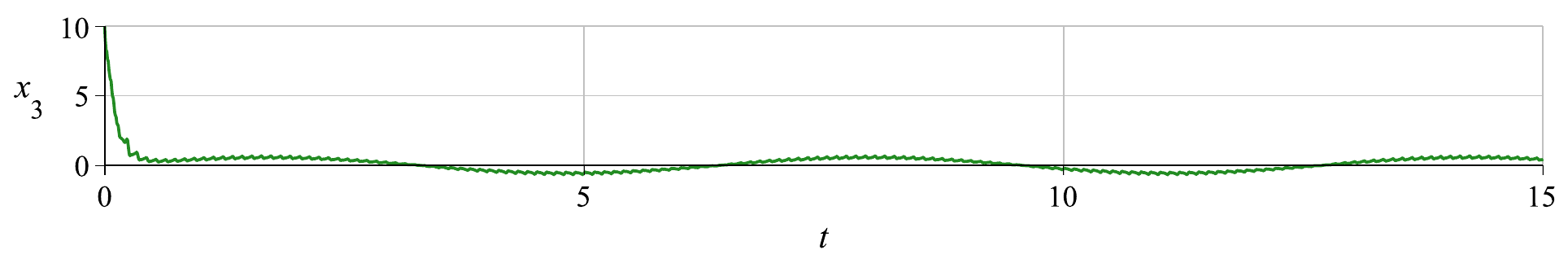}
\includegraphics[width=1\linewidth]{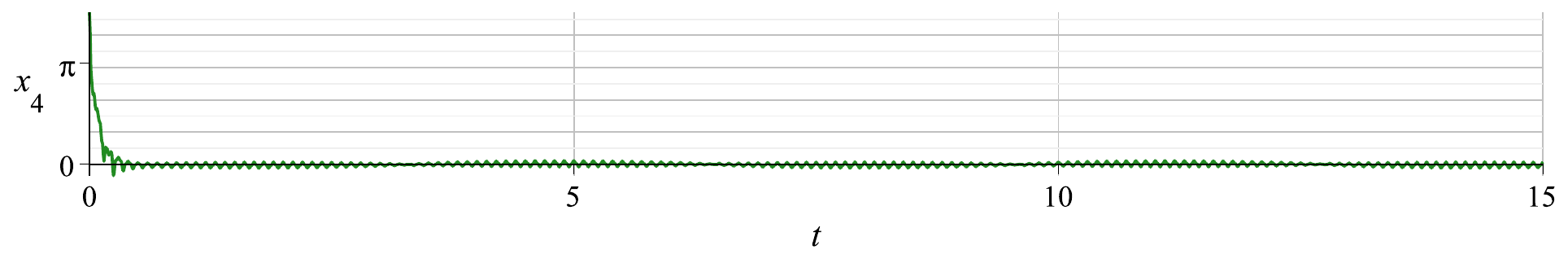}
\includegraphics[width=1.03\linewidth]{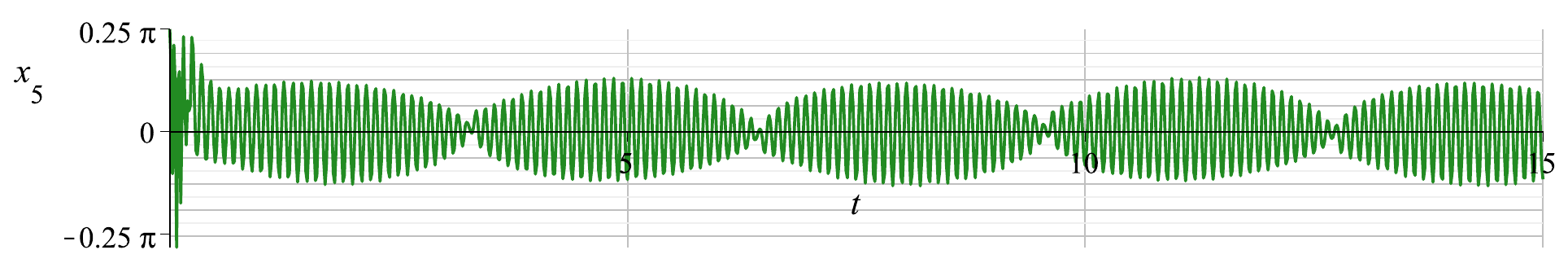}
\includegraphics[width=1\linewidth]{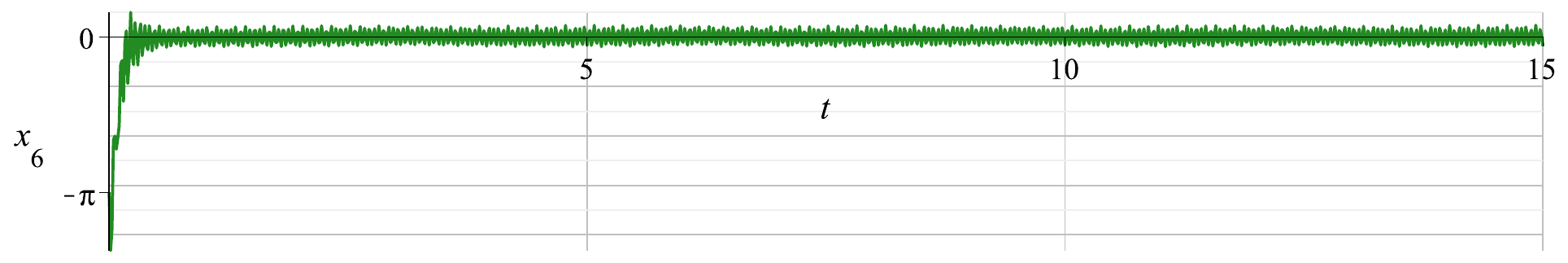}
a)
\end{center}
 \end{minipage}\hfill
 \begin{minipage}{0.49\linewidth}
\begin{center}
\includegraphics[width=1\linewidth]{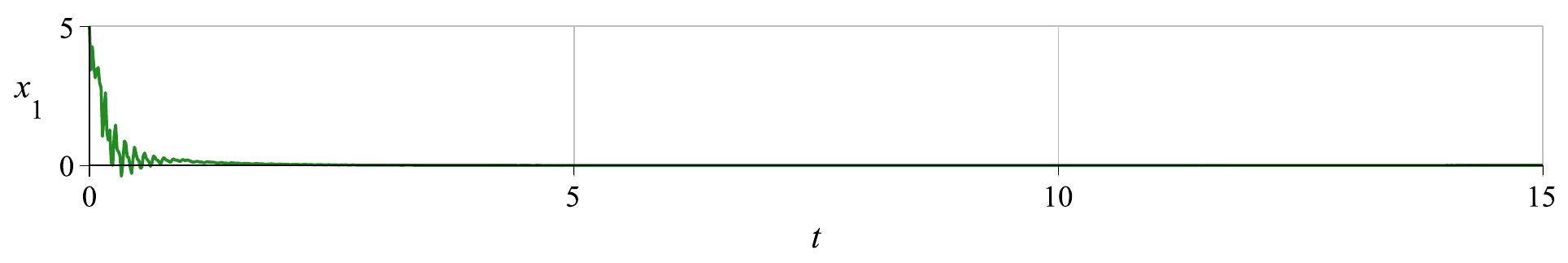}
\includegraphics[width=1\linewidth]{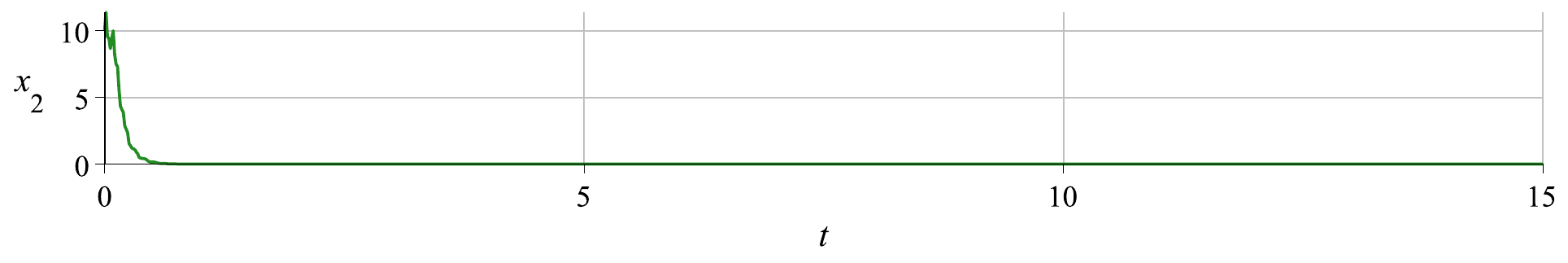}
\includegraphics[width=1\linewidth]{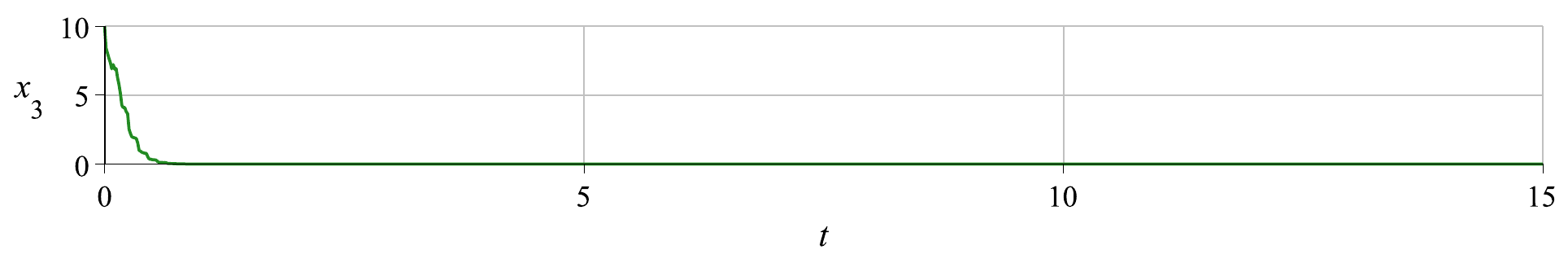}
\includegraphics[width=1\linewidth]{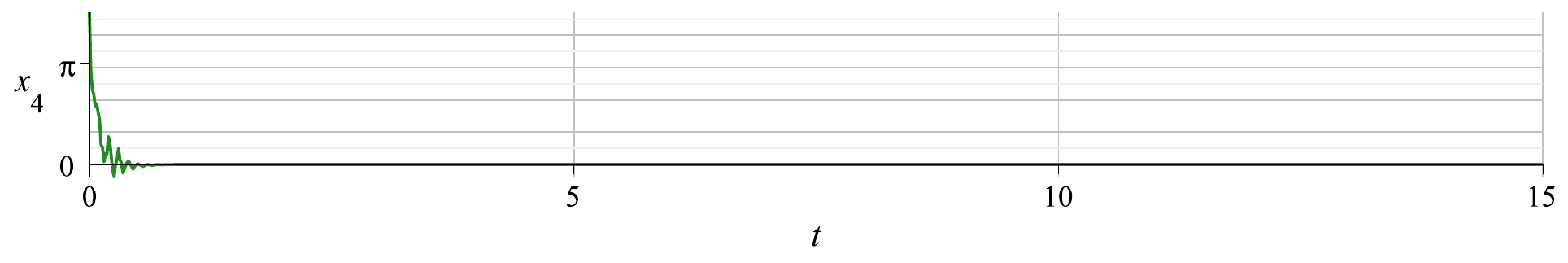}
\includegraphics[width=1.03\linewidth]{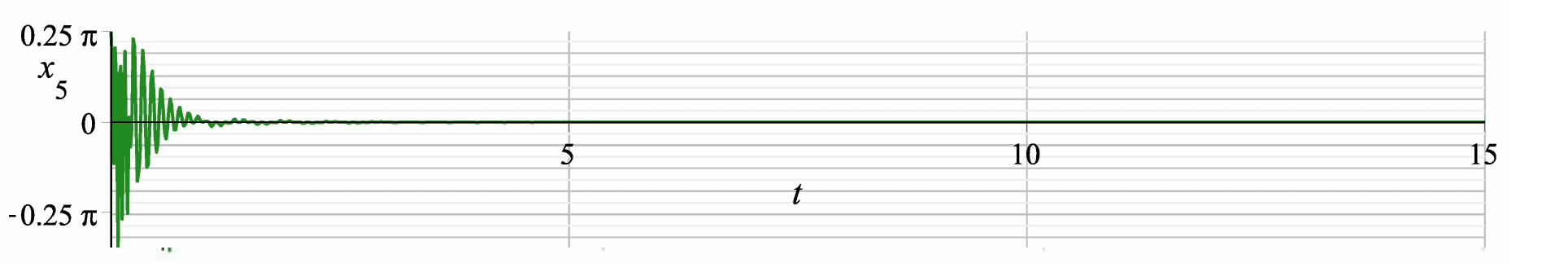}
\includegraphics[width=1\linewidth]{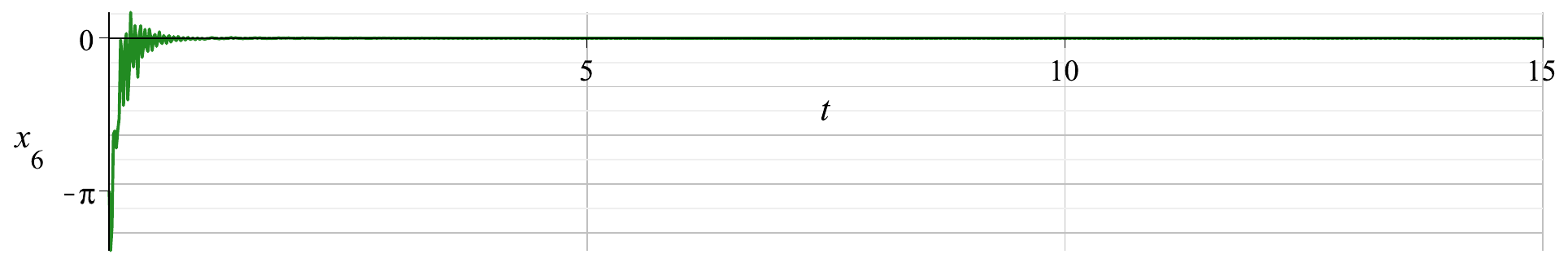}
b)
\end{center}
 \end{minipage}
\caption{Time-plots of the trajectories of system~\eqref{underdr} with controls~\eqref{cont_ship}.}
\end{figure*}
Consider the equations of motion for an autonomous 3D
underwater vehicle studied, e.g., in~\cite{Bara}, and assume that the motion of the vehicle is also affected by external disturbances:
\begin{equation}\label{underdr}
\begin{aligned}
\dot x_1=\sum_{k=1}^4f_k(x)u_k+g(t),
\end{aligned}
 \end{equation}
 where   $(x_1, x_2, x_3)$  are the coordinates of the center of mass, $(x_4$, $x_5$, $x_6)$ describe the vehicle orientation (Euler angles),  $u_1$ is the translational velocity  along the $Ox_1$ axis,  $(u_2,u_3,u_4)$ are the angular velocity components, and the vector fields of the unperturbed system are
$$
 \begin{aligned}
&f_1(x)=(\cos x_5\cos x_6, \cos x_5\sin x_6,{-}\sin x_5,0,0,0)^\top, \\
&f_2(x) =(0,0,0,1,0,0)^\top,\\
 & f_3(x){=}(0,0,0,\sin x_4{\rm tg}\, x_5,\cos x_4,\sin x_4\sec x_5)^\top,\\
 &f_4(x){=}(0,0,0,\cos x_4{\rm tg}\, x_5,{-}\sin x_4,\cos x_4\sec x_5)^\top.
 \end{aligned}
$$
The drift term in~\eqref{underdr} accounts for the external disturbances caused by waves and ocean currents, and we choose the following form for $g(t)$:
 $$
g(t)=(0,d,a\sin(\omega t+b),0,0,0)^\top,
$$
where $a,b,d,\omega$ are some positive constants.
The rank condition~\eqref{rank} is satisfied in the domain $D =\{x \in \mathbb R^6 \,|\,-\tfrac{\pi}{2}<x_5<\tfrac{\pi}{2}\}$ with $S_1=\{1,2,3,4\}$, $S_2=\{(1,3),(1,4)\}$, $S_3=\emptyset$. Then the matrix~\eqref{fmatrix}  takes the form
$$
\mathcal F(x)=\left(
f_1(x),\, f_2(x),\, f_3(x),\, f_4(x),\, [f_1,f_3](x),\, [f_1,f_4](x)
              \right),
$$
and we may write controls~\eqref{cont} as $u_k=h_k^\varepsilon(t,x)$:
 \begin{align}
    h_1^\varepsilon(t,x)=&a_1(x)+2{\ \rm sign}(a_{13}(x))\sqrt{\tfrac{\pi|a_{13}(x)|}{\varepsilon}}\cos\tfrac{2\pi k_{13}t}{\varepsilon}\nonumber\\
   &+2{\ \rm sign}(a_{14}(x))\sqrt{\tfrac{\pi|a_{14}(x)|}{\varepsilon}}\cos\tfrac{2\pi k_{14}t}{\varepsilon},\nonumber\\
   h_2^\varepsilon(t,x)=&a_2(x),\label{cont_ship}\\
    h_3^\varepsilon(t,x)=&a_3(x)+2\sqrt{\tfrac{\pi|a_{13}(x)|}{\varepsilon}}\sin\tfrac{2\pi k_{13} t}{\varepsilon},\nonumber
    \end{align}
     \begin{align}
    h_4^\varepsilon(t,x)=&a_4(x)+2\sqrt{\tfrac{\pi|a_{14}(x)|}{\varepsilon}}\sin\tfrac{2\pi k_{14} t}{\varepsilon},\nonumber
 \end{align}
with
$
a(x)=\left(a_1(x),\,a_2(x),\,a_{13}(x),\,a_{14}(x)\right)^\top{=}-\gamma \mathcal F^{-1}(x)x.
$
The behavior of  system~\eqref{underdr} with controls~\eqref{cont_ship} is illustrated in Fig.~1a). For numerical simulations, we take
$g(t)=(0,2,5\sin t,0,0,0)^\top,$ $x^0=\Big(5,10,10,\frac{3\pi}{2},\frac{\pi}{4},-\pi\Big)^\top$, $\varepsilon=0.1$, $\gamma=10$, $k_{13}=1$, $k_{14}=2$.
To illustrate Corollary~1, assume that the drift is described by $g(t,x)=(0,x_1^3(t),x_2^3(t)\sin t,0,0,0)^\top.$ As it is shown in Fig.1b), the trajectories of system~\eqref{underdr} tend asymptotically to zero in this case.
\subsection{Front-wheel drive car}
As an example of a nonholonomic system satisfying condition~\eqref{rank} with the second-order Lie brackets, consider a kinematic model of the front-wheel drive car (see, e.g.,~\cite{Luca95}):
$$
\begin{aligned}
\dot x_1=\sum_{k=1}^2f_k(x)u_k,
\end{aligned}
$$
 where $(x_1,x_2)$ are  the
Cartesian coordinates of the rear axle center, the angle $x_3$ defines the car orientation with respect to the $x_1$-axis, $x_4$ is the steering angle, $u_1,u_2$ denote the driving and the steering velocity input, respectively; thus the vector fields of the system are given by
$$
 \begin{aligned}
&f_1(x)=(\cos x_3\cos x_4,\sin x_3\cos x_4,\sin x_4,0)^\top, \\
&f_2(x) =(0,0,0,1)^\top.\\
 \end{aligned}
$$
It can be verified that the rank condition~\eqref{rank} is satisfied with $S_1=\{1,2\}$, $S_2=\{(1,2)\}$, $S_3=\{(1,2,1)\}$, so that the matrix
$$
\mathcal F(x){=}\left(
f_1(x), f_2(x),  [f_1,f_2](x), \big[[f_1,f_2],f_1\big](x)
              \right)
$$
is nonsingular in $\mathbb R^4$. If the control input acts with an error, i.e. $u_k=h_k^\varepsilon(t,x)+n_k(t,x)$, where $n_k(t,x)$ are some disturbances, then the system equations can be interpreted as the system with drift:
\begin{equation}\label{ex_car}
  \begin{aligned}
\dot x_1=\sum_{k=1}^2f_k(x)u_k+g(t,x),
\end{aligned}
\end{equation}
where $g(t,x)=\sum\limits_{k=1}^2f_k(x)n_k(x,t)$. According to the proposed design procedure, we take
 controls of the form~\eqref{cont}:
 \begin{align}
h_1^\varepsilon(t,x)=&a_1(x)+2{\ \rm sign}(a_{12}(x))\sqrt{\tfrac{\pi|a_{12}(x)|}{\varepsilon}}\cos\tfrac{2\pi k_{12}t}{\varepsilon}\nonumber\\
  &+2\sqrt[3]{\tfrac{2\pi^2(k_{2121}-k_{1121})a_{121}(x)}{\varepsilon^2}}\cos\tfrac{2\pi k_{1121}t}{\varepsilon}\nonumber\\
  &\times\Big(1+\tfrac{2\pi k_{2121}t}{\varepsilon}\sin\Big),\label{cont_car}\\
h_2^\varepsilon(t,x)=&a_2(x)+2\sqrt{\tfrac{\pi|a_{12}(x)|}{\varepsilon}}\sin\tfrac{2\pi k_{12}t}{\varepsilon},\nonumber\\
  &+2\sqrt[3]{\tfrac{2\pi^2(k_{2121}-k_{1121})a_{121}(x)}{\varepsilon^2}}\sin\tfrac{2\pi k_{2121}t}{\varepsilon}\nonumber
 \end{align}
with
$$
\begin{aligned}
a(x)&=\left(a_1(x),\,a_2(x),\,a_{12}(x),\,a_{121}(x)\right)^\top=-\gamma \mathcal F^{-1}(x)x.
\end{aligned}
$$
For the numerical simulation, we take $n_1(t,x)=2\cos 10\pi t$, $n_2(t,x)=\sin 20\pi t$, $x^0=\Big(5,3,-\frac{\pi}{2},\frac{\pi}{4}\Big)^\top$, $\varepsilon=0.5$, $\gamma=15$, $k_{12}=7$, $k_{1121}=3$, $k_{2121}=1$. The corresponding plots are depicted in Fig.~2.
\begin{figure}[tpt]
\includegraphics[width=1\linewidth]{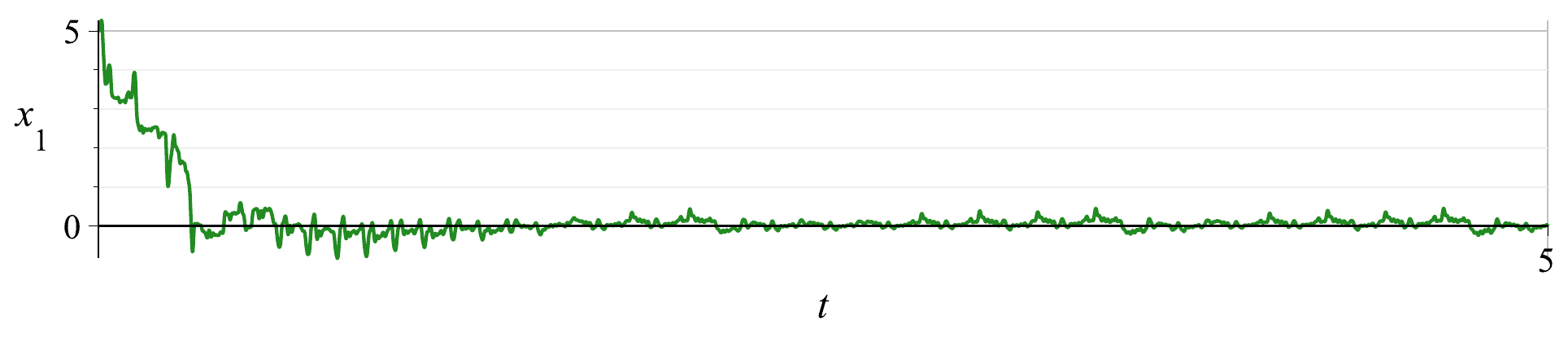}
\includegraphics[width=1\linewidth]{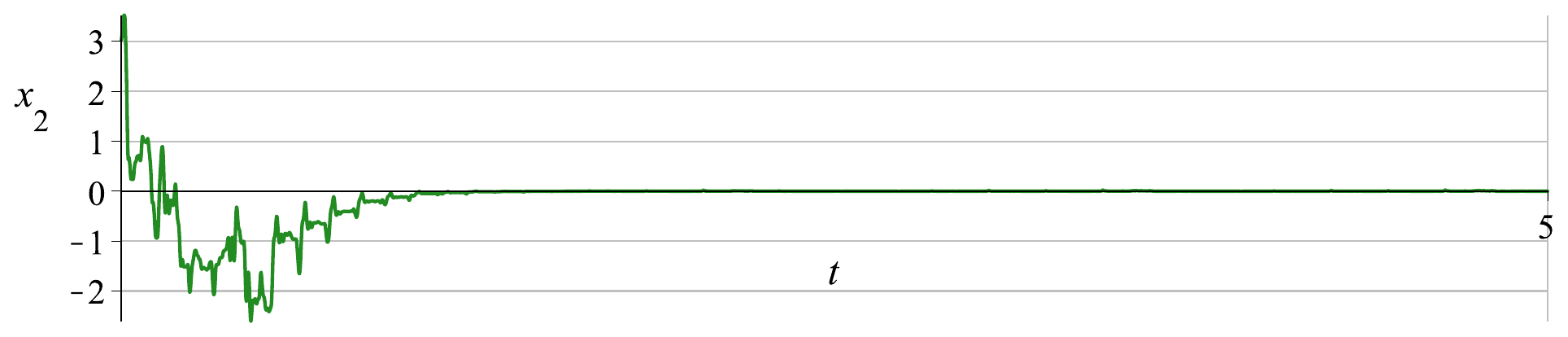}
\includegraphics[width=1\linewidth]{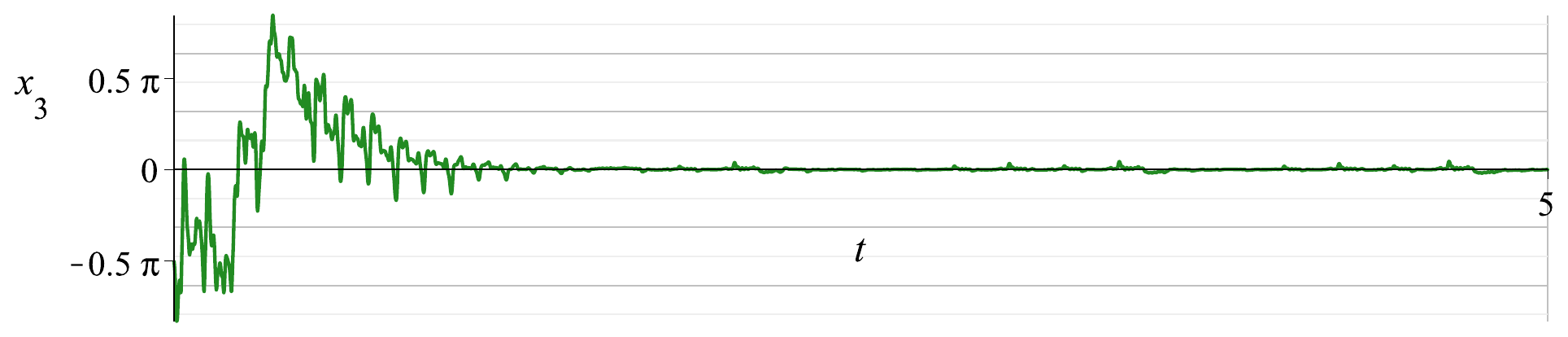}
\includegraphics[width=1\linewidth]{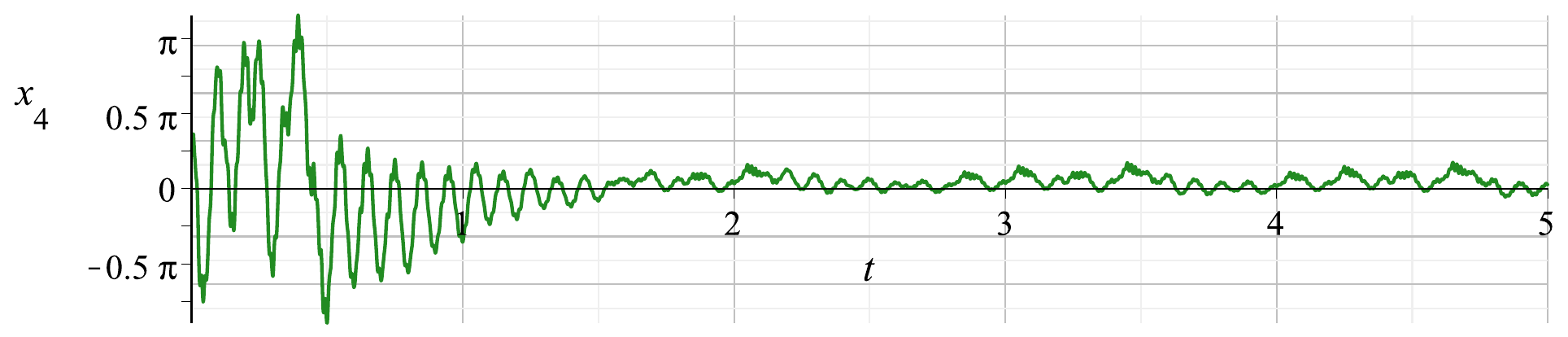}
\caption{Time-plots of the trajectories of system~\eqref{ex_car} with controls~\eqref{cont_car}.}
\end{figure}
\section{Conclusions}
We have considered a class of nonholonomic systems with  time-varying drift term satisfying certain boundedness assumptions.
Extending the approach of~\cite{ZGB16,GZ18}, we have obtained a family of time-periodic control functions with rather simple  formulas for state-dependent coefficients. It should be emphasized that the considered systems with vanishing controls, in general, do  not admit the trivial equilibrium. It is also crucial that the exponential decay estimates have been derived without assuming that the drift can be compensated by a linear combination of control vector fields.


\begin{thebibliography}{10}
\bibitem{Ag12}
Aguilar, C.O. (2012).
\newblock Local controllability of control-affine systems with quadractic drift
  and constant control-input vector fields.
\newblock In \emph{Proc. 51st IEEE Conference on Decision and Control},
  1877--1882.

\bibitem{Ast94}
Astolfi, A. (1994).
\newblock On the stabilization of nonholonomic systems.
\newblock In \emph{Proc. 33rd IEEE Conference on Decision and Control},
  volume~4, 3481--3486.

\bibitem{Bara}
Barraquand, J. and Latombe, J.C. (1989).
\newblock On nonholonomic mobile robots and optimal maneuvering.
\newblock In \emph{Proc. IEEE International Symposium on Intelligent Control},
  340--347.

\bibitem{Bu00}
Bullo, F., Leonard, N.E., and Lewis, A.D. (2000).
\newblock Controllability and motion algorithms for underactuated Lagrangian
  systems on Lie groups.
\newblock \emph{IEEE Transactions on Automatic Control}, 45(8), 1437--1454.

\bibitem{Clar97}
Clarke, F.H., Ledyaev, Y.S., Sontag, E.D., and Subbotin, A.I. (1997).
\newblock Asymptotic controllability implies feedback stabilization.
\newblock \emph{IEEE Transactions on Automatic Control}, 42(10), 1394--1407.

\bibitem{Cor92}
Coron, J.M. (1992).
\newblock Global asymptotic stabilization for controllable systems without
  drift.
\newblock \emph{Mathematics of Control, Signals, and Systems}, 5(3), 295--312.

\bibitem{Luca95}
De~Luca, A. and Oriolo, G. (1995).
\newblock Modelling and control of nonholonomic mechanical systems.
\newblock In \emph{Kinematics and Dynamics of Multi-Body Systems}, 277--342.
  Springer.

\bibitem{Thi00}
Floquet, T., Barbot, J.P., and Perruquetti, W. (2000).
\newblock One-chained form and sliding mode stabilization for a nonholonomic
  perturbed system.
\newblock In \emph{Proc. 2000 American Control Conference}, volume~5,
  3264--3268.

\bibitem{Gao11}
Gao, F., Yuan, F., Yao, H., and Mu, X. (2011).
\newblock Adaptive stabilization of high order nonholonomic systems with strong
  nonlinear drifts.
\newblock \emph{Applied Mathematical Modelling}, 35(9), 4222--4233.

\bibitem{Godv99}
Godhavn, J.M., Balluchi, A., Crawford, L.S., and Sastry, S.S. (1999).
\newblock Steering of a class of nonholonomic systems with drift terms.
\newblock \emph{Automatica}, 35(5), 837--847.

\bibitem{GZ18}
Grushkovskaya, V. and Zuyev, A. (2018).
\newblock Obstacle avoidance problem for second degree nonholonomic systems.
\newblock In \emph{Proc. 57th IEEE Conf. on Decision and Control}, 1500--1505.

\bibitem{GZE18}
Grushkovskaya, V., Zuyev, A., and Ebenbauer, C. (2018).
\newblock On a class of generating vector fields for the extremum seeking
  problem: Lie bracket approximation and stability properties.
\newblock \emph{Automatica}, 94, 151--160.

\bibitem{Her80}
Hermes, H. (1980).
\newblock On the synthesis of a stabilizing feedback control via Lie algebraic
  methods.
\newblock \emph{SIAM Journal on Control and Optimization}, 18(4), 352--361.

\bibitem{Jean15}
Jean, F. and Prandi, D. (2015).
\newblock Complexity of control-affine motion planning.
\newblock \emph{SIAM Journal on Control and Optimization}, 53(2), 816--844.

\bibitem{Kolman95}
Kolmanovsky, I. and McClamroch, N.H. (1995).
\newblock Developments in nonholonomic control problems.
\newblock \emph{IEEE Control Systems}, 15(6), 20--36.

\bibitem{MCl98}
M'Closkey, R. and Morin, P. (1998).
\newblock Time-varying homogeneous feedback: design tools for the exponential
  stabilization of systems with drift.
\newblock \emph{International Journal of Control}, 71(5), 837--869.

\bibitem{Mich03}
Michalska, H. and Torres-Torriti, M. (2003).
\newblock A geometric approach to feedback stabilization of nonlinear systems
  with drift.
\newblock \emph{Systems \& Control Letters}, 50(4), 303--318.

\bibitem{Mor99}
Morin, P., Pomet, J.B., and Samson, C. (1999).
\newblock Design of homogeneous time-varying stabilizing control laws for
  driftless controllable systems via oscillatory approximation of Lie brackets
  in closed loop.
\newblock \emph{SIAM Journal on Control and Optimization}, 38(1), 22--49.

\bibitem{Pom99}
Pomet, J.B. (1999).
\newblock On the curves that may be approached by trajectories of a smooth
  control affine system1.
\newblock \emph{Systems \& Control Letters}, 36(2), 143--149.

\bibitem{Rey99}
Reyhanoglu, M., van~der Schaft, A., McClamroch, N.H., and Kolmanovsky, I.
  (1999).
\newblock Dynamics and control of a class of underactuated mechanical systems.
\newblock \emph{IEEE Transactions on Automatic Control}, 44(9), 1663--1671.

\bibitem{Wang04}
Wang, Z., Ge, S., and Lee, T. (2004).
\newblock Adaptive stabilization of nonholonomic systems with unknown virtual
  control coefficients and nonlinear drifts.
\newblock In \emph{Proc. of the 2004 IEEE International Symposium on
  Intelligent Control}, 281--286.

\bibitem{YaYa10}
Yang, L. and Yang, J. (2010).
\newblock Stabilization for a class of nonholonomic perturbed systems via
  robust adaptive sliding mode control.
\newblock In \emph{Proc. American Control Conference 2010}, 1178--1183.

\bibitem{Zhao13}
Zhao, Y. and Wu, Y. (2013).
\newblock Control of nonholonomic systems with nonlinear unmeasured dynamics by
  output feedback.
\newblock \emph{Journal of Control Theory and Applications}, 11(3), 504--512.

\bibitem{ZuSIAM}
Zuyev, A. (2016).
\newblock Exponential stabilization of nonholonomic systems by means of
  oscillating controls.
\newblock \emph{SIAM Journal on Control and Optimization}, 54(3), 1678--1696.

\bibitem{ZG17}
Zuyev, A. and Grushkovskaya, V. (2017).
\newblock Motion planning for control-affine systems satisfying low-order
  controllability conditions.
\newblock \emph{International Journal of Control}, 90, 2517--2537.

\bibitem{ZG19}
Zuyev, A. and Grushkovskaya, V. (2019).
\newblock On stabilization of nonlinear systems with drift by time-varying
  feedback laws.
\newblock To appear in \emph{Proc. 12th Int. Workshop on Robot Motion and
  Control}, arXiv:1703.02348.

\bibitem{ZGB16}
Zuyev, A., Grushkovskaya, V., and Benner, P. (2016).
\newblock Time-varying stabilization of a class of driftless systems satisfying
  second-order controllability conditions.
\newblock In \emph{Proc. 15th European Control Conference}, 1678--1696.
\end{thebibliography}
\end{document}